\documentclass{nmeauth}

\usepackage{algorithmic}
\usepackage{tabularx,subfigure,caption}
\usepackage{tikz}
\usepackage[active]{srcltx}
\usepackage{type1cm}
\usepackage{comment}
\usepackage{color}
\usepackage{url}
\usepackage{undertilde}

\newcommand{\stiff}{\ensuremath{K}}
\newcommand{\dep}{\ensuremath{u}}

\newcommand{\force}{\ensuremath{f}}
\newcommand{\lam}{\ensuremath{\lambda_b}}

\newcommand{\glam}{\ensuremath{\lambda}}
\newcommand{\drhs}{\ensuremath{d}}
\newcommand{\de}{\ensuremath{e}}
\newcommand{\dgschur}{\ensuremath{F}}
\newcommand{\dG}{\ensuremath{G}}

\newcommand{\dalpha}{\ensuremath{g}}
\newcommand{\domain}{\ensuremath{\Omega}}
\newcommand{\s}{\ensuremath{^{(s)}}}
\newcommand{\g}{\ensuremath{^{\domain}}}
\newcommand{\n}{\ensuremath{^{()}}}
\newcommand{\assem}{\ensuremath{A}}
\newcommand{\dassem}{\ensuremath{A}}

\newcommand{\krylov}{\ensuremath{\mathcal{K}}}

\newcommand{\schur}{\ensuremath{S}}
\newcommand{\kernel}{\ensuremath{N}}

\newcommand{\projrig}{\ensuremath{Q}}
\newcommand{\scal}{\ensuremath{\Delta}}
\newcommand{\range}{\ensuremath{\operatorname{range}}}
\newcommand{\spann}{\ensuremath{\operatorname{span}}}
\newcommand{\diag}{\ensuremath{\operatorname{diag}}}

\newcommand{\uu}[1]{\ensuremath{\underline{#1}}}
\newcommand{\uuu}[1]{\ensuremath{\underline{\underline{#1}}}}
\newcommand{\muv}[1]{\ensuremath{\utilde{#1}}}
\newcommand{\trace}{\ensuremath{t}}

\newcommand{\id}{\ensuremath{\operatorname{I}}}

\begin{document}
\NME{1}{6}{00}{28}{00}

\runningheads{Gosselet, Rixen and Rey}{Structures containing repeated patterns}

\title{A domain decomposition strategy to efficiently solve structures containing repeated patterns}

\author{P. Gosselet\affil{1}\corrauth, D. J. Rixen\affil{2}, C. Rey\affil{1}}

\address{\affilnum{1} LMT Cachan, ENS Cachan/CNRS/UPMC/PRES UniverSud Paris,\\ 61 av. du pr\'esident Wilson, F-94230 CACHAN, FRANCE\\
\affilnum{2} Delft University of Technology, Faculty 3mE, Department of Precision and Microsystems Engineering,
Engineering Dynamics, Mekelweg 2, 2628 CD DELFT, THE NETHERLANDS
}

\corraddr{gosselet@lmt.ens-cachan.fr}



\received{}
\revised{}
\noaccepted{}

\begin{abstract}
This paper presents a strategy for the computation of structures with repeated patterns based on
domain decomposition and block Krylov solvers. It can be seen as a special variant of the FETI method.
We propose using the presence of repeated domains in the problem to
compute the solution by minimizing the interface error on several directions simultaneously.
The method not only drastically decreases the size of the problems to solve but also accelerates the convergence
of interface problem for nearly no additional computational cost and minimizes expensive memory accesses.
The numerical performances are illustrated on some thermal and elastic academic problems.
\end{abstract}


\section{INTRODUCTION}
Industrial structures often exhibit repeated patterns. Indeed the use of repetitions simplifies the design and
enables to reduce costs in assembled structures because repeated components can be mass-produced;
in addition symmetries enable to equilibrate inertia and to balance loadings. The existence of repetitions
is naturally used when defining the CAD model of the structure, using ``copy-paste'' and rigid transformations.
Unfortunately, taking into account the repetitions during the numerical simulation of the structure is not as an easy task.
Strategies based on Fourier expansions \cite{RIXEN.2005.1} or homogenization techniques \cite{Bergman.1985.1} are
only valid in a limited range of application, most of them require either true geometrical periodicity,
repetition of one single pattern or large number of repetitions;
their application to non-periodic geometries, loadings and boundary conditions is not straightforward.

\begin{figure}[ht]
\centering\includegraphics[width=5cm]{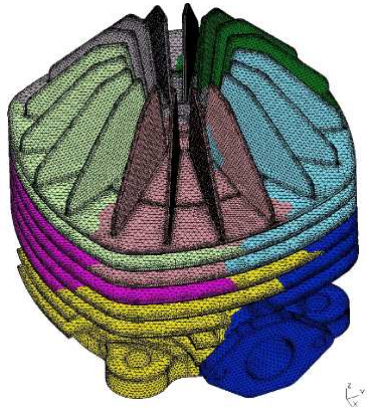}
\caption{EDF carter (courtesy of O. Boiteau)}\label{fig:edf}
\end{figure}

The crankcase shown in figure \ref{fig:edf} and typically used in power stations of EDF (Electricit\'e de France) is a good example of the structures we dedicate our method to:
it is composed of repeated cooling winglets and disks, and a complex clamping system.
The winglets and disks are repeated patterns and we will assume that the mesh of every occurrence of a given pattern
has the same discretization (as is often the case in practice to simplify the mesh generation). The remainder of the structure is not a repeated pattern.
The winglets and the disk
have different orientations and thus the structure in its entirety cannot be considered as quasi-cyclic.
The classical method to simulate the carter would be to obtain a complete mesh of the structure and to solve
this large system (possibly using a domain decomposition algorithm with automatic substructuring software).

\begin{figure}[ht]
\centering\includegraphics[width=7cm]{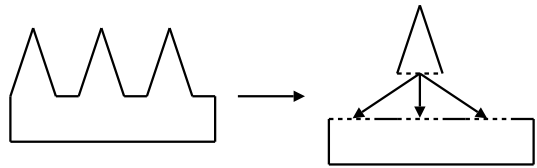}
\caption{Schematic view of structures with occurrences of a pattern}\label{fig:simple}
\end{figure}

In the present work we propose a numerical strategy based on domain decomposition methods \cite{FARHAT:1994:ADV,LETALLEC:1994:DDM,FARHAT:2001:FETI_DP} (see \cite{GOSSELET.2007.1} for a review)
and block-Krylov solvers \cite{SAAD:2000:IMS,ALIAGA.1999.1,GUTKNECHT.2005.1}.
The main idea behind the strategy developed here is to redefine the structure as a collection of occurrences
of patterns as described in figure \ref{fig:simple}: this structure is made of two patterns, one of which occurs
three times. Our method exhibits the following properties compared to classical domain decomposition methods:
\begin{itemize}
\item only the patterns need to be meshed (not their occurrences),
\item the computational operations are realized more efficiently,
\item the iteration schemes takes advantage of the numerical information associated to all occurrences to converge faster.
\end{itemize}

The paper is organized as follow. First, basics on domain decomposition and Krylov solvers are briefly recalled.
Then the method is presented on the most simple case: the linear thermal analysis of a fully periodic structure (one repeated pattern). The method is then generalized to non-periodic structures, to elastic problems and to decompositions involving floating substructures. Eventually, conclusion and prospects are provided.


\section{Basics on domain decomposition methods and Krylov solvers}
The method we propose is a variant of the most classical non-overlapping domain decomposition methods
\cite{GOSSELET.2007.1}. For the sake of simplicity, we restrict our presentation to the dual domain
decomposition method (a.k.a. FETI method \cite{FARHAT:1994:ADV,GOSSELET:2003:IEI}) in the context of the finite element analysis of a linear elasticity problem.

Let us denote \stiff\g\, the stiffness matrix, \dep\g\, the displacement field and \force\g\, the generalized forces
set on domain \domain. The so called ``global''(or ``assembled'') problem to be solved is
\begin{equation}\label{eq:equilibre}
\stiff\g \dep\g = \force\g
\end{equation}
We consider a partition of \domain\, into $n$ non-overlapping subdomains denoted \domain\s,
and introduce interface internal forces \lam\s\, imposed on subdomain \domain\s\, by its neighbors.
We call \trace\s\, the trace operator which extracts boundary degrees of freedom from subdomain \domain\s.
Equilibrium of subdomains, equilibrium of interfaces and connectivity of submeshes read:
\begin{equation}\label{eq:equilibre2}
\begin{aligned}
&\stiff\s \dep\s = \force\s - {\trace\s}^T \lam\s \qquad \text{for }s=1\ldots n\\
&\sum_{s=1}^n \dassem\s \trace\s \dep\s = 0 \\
&\lam\s ={\dassem\s}^T \glam \qquad \text{for }s=1\ldots n
\end{aligned}
\end{equation}
where \dassem\s\, is the signed boolean assembly operator which connects pairwise degrees of freedom on the interface. We have chosen to express the reactions on the interface of subdomains from one unique interface stress field \glam\, insuring automatically the equilibrium of the interface reactions. The second equation corresponds to the equality of interface displacements. In the following subscript $b$ stands for boundary and subscript $i$ for internal, so that $\trace\s \dep\s=\dep_{b}\s$.

In case the Dirichlet boundary conditions of a subdomain are not sufficient to fix it, we call the subdomain a ``floating subdomain''; the kernel of the associated stiffness matrix $\stiff\s$ is denoted by $\kernel\s$, a pseudo-inverse of $\stiff\s$ is denoted by ${\stiff\s}^+$, it verifies $\stiff\s{\stiff\s}^+ y = y$ for any $y$ in $\range(\stiff\s)$.
Eliminating the displacements $\dep\s$ from the formulation (\ref{eq:equilibre2}), one finds the dual interface problem
as appearing in the Finite Element Tearing and Interconnecting (FETI) method \cite{FARHAT:1994:ADV}:
\begin{equation}\label{eq:DUALequilibre}
  \begin{pmatrix} \dgschur & \dG \\ \dG^T & 0 \end{pmatrix} \begin{pmatrix} \glam \\ \dalpha \end{pmatrix} = \begin{pmatrix} \drhs \\ \de \end{pmatrix}
\end{equation}
where we introduced the following definitions:
\begin{equation*}
 \begin{aligned}
  \dgschur = \sum_s \dassem\s \trace\s {\stiff\s}^+{\trace\s}^T{\dassem\s}^T & \qquad & \dG=\begin{pmatrix}\ldots\dassem\s\trace\s\kernel\s\ldots\end{pmatrix}\\
\drhs=\sum_s \dassem\s\trace\s{\stiff\s}^+\force\s & & \de^T=\begin{pmatrix}\ldots{\force\s}^T\kernel\s\ldots\end{pmatrix}
 \end{aligned}
\end{equation*}
We will also denote local Schur complement by $\schur\s=\stiff\s_{bb}-\stiff\s_{bi}{\stiff\s_{ii}}^{-1}\stiff\s_{ib}$.
It can easily be proven that ${\schur\s}^+=\trace\s {\stiff\s}^+{\trace\s}^T$ and thus
$\dgschur = \sum_s \dassem\s {\schur\s}^+ {\dassem\s}^T$.
In order to take advantage of the additive structure of matrix \dgschur, a Krylov iterative solver is chosen to solve the
dual interface problem (\ref{eq:DUALequilibre}).
Our strategy is independent of the Krylov solver. For simplicity reasons we choose, as in the basic FETI methods, to
use the conjugate gradient. We hence suppose that \dgschur\, is symmetric positive definite (properties which
are directly inherited from matrix \stiff\g).

One key point of coupling domain decomposition methods and
Krylov iterative solvers is the preconditioner $\tilde{\dgschur}^{-1}$.
 Various versions have been proposed for the FETI method:
 the so called ``Dirichlet'', ``lumped'', ``superlumped'' and ``identity`` preconditioners
 (in increasing order of efficiency and computational burden):
\begin{align}
\tilde{\dgschur}^{-1}_D &= \sum_s \scal\s{\dassem}\s \schur\s {\dassem\s}^T\scal\s \qquad&\qquad
\tilde{\dgschur}^{-1}_L &= \sum_s \scal\s{\dassem}\s\stiff\s_{bb}{\dassem\s}^T\scal\s \\
\tilde{\dgschur}^{-1}_{SL} &= \sum_s \scal\s{\dassem}\s\diag({\stiff}\s_{bb}){\dassem\s}^T\scal\s\qquad&\qquad\
\tilde{\dgschur}^{-1}_{I} &= \id
\end{align}
where $\scal\s$  is a diagonal scaling operator ($\sum_s \scal\s = I$).

Another keypoint in solving the dual interface problem is the coarse problem imposing at every
iteration the self-equilibrium constraint of the interface forces, the second set of equations
in (\ref{eq:DUALequilibre}). The coarse problem is solved using a projector $\projrig$ and initialization $\glam_0$,
defined by
\begin{align}
\projrig &= I - \tilde{\dgschur}^{-1}\dG\left(\dG^T \tilde{\dgschur}^{-1} \dG\right)^{-1}\dG^T  \\
\glam_0 &=  \tilde{\dgschur}^{-1}\dG\left(\dG^T \tilde{\dgschur}^{-1} \dG\right)^{-1}\de
\end{align}
where $\tilde{\dgschur}^{-1}$ can be any of the preconditioners introduced above.

The standard FETI solver is summarized in the left column of table \ref{tab:alggc}.
\begin{table}[ht]\caption{Dual Schur complement with conjugate gradient}\label{tab:alggc}
\centering
\begin{tabular}{|l||m{6cm}|m{6cm}|}\hline
  & \textbf{Classical version} & \textbf{Multivector version} \\ \hline
1 & \multicolumn{2}{c|}{Initialize with arbitrary $\glam_{00}$}\\ \hline
2 & \multicolumn{2}{c|}{Compute $\glam_0 = \tilde{\dgschur}^{-1}\dG (\dG^T \tilde{\dgschur}^{-1} \dG)^{-1} \de + Q \glam_{00}$ } \\ \hline
3 & \multicolumn{2}{c|}{Compute $r_0=Q^T (\drhs -\dgschur \glam_0)$} \\ \hline
4 & $z_0 = Q\tilde{\dgschur}^{-1}r_0$ set $w_0=z_0$ & $z_0 = Q\tilde{\dgschur}^{-1}r_0$ set $w_0=z_0$ obtain $\muv{w}_0$ through permutations \\ \hline
5 & \multicolumn{2}{c|}{for $j=0,\ldots,m$}\\ \hline
5.1 & $p_j = Q^T\dgschur  w_j $ & $\muv{p}_j = Q^T\dgschur  \muv{w}_j $ \\ \hline
5.2 & $\alpha_j=(w_j^T r_j)/(w_j^T p_j)$ & $\uu{\alpha}_j=({\muv{w}_j}^T{\muv{p}_j})^{-1}({\muv{w}_j}^Tr_j)$\\ \hline
5.3 & $\glam_{j+1}=\glam_j+\alpha_j w_j$ & $\glam_{j+1}=\glam_j+\muv{w}_j\uu{\alpha}_j$\\ \hline
5.4 & $r_{j+1}=r_j-\alpha_j p_j$ & $r_{j+1}=r_j-\muv{p}_j\uu{\alpha}_j$\\ \hline
5.5 & $z_{j+1}=Q\tilde{\dgschur}^{-1}r_{j+1}$&$z_{j+1}=Q\tilde{\dgschur}^{-1}r_{j+1}$ obtain $\muv{z}_{j+1}$ through permutations\\ \hline
5.6 & $\beta_j=-(p_j^T z_{j+1})/(w_j^T p_j)$ & $\uuu{\beta}_j=-({\muv{w}_j}^T\muv{p}_j)^{-1}(\muv{p}_j^T{\muv{z}_{j+1}})$ \\ \hline
5.7 &$w_{j+1}=z_{j+1}+\beta_j w_j$ &  $\muv{w}_{j+1}=\muv{z}_{j+1}+\muv{w}_j\uuu{\beta}_j $ \\ \hline
\end{tabular}
\end{table}

Conjugate gradient consists in building a $\dgschur$-orthogonal basis $(w_i)_{0\leqslant i\leqslant m}$ of Krylov subspace $\krylov_m\left(z_0,Q\tilde{\dgschur}^{-1}Q^T\dgschur \right)=\spann\left(z_0, \ldots, \left(Q\tilde{\dgschur}^{-1}Q^T\dgschur\right)^{m-1} z_0\right)$ and finding approximation $\glam_m\in \glam_0+\krylov_m\left(z_0,Q\tilde{\dgschur}^{-1}Q^T\dgschur \right)$ so that at each iteration the error $\|\glam-\glam_m\|_\dgschur$ is minimized (which is equivalent to making the residual $r_m=(d-\dgschur\glam_m)$ orthogonal to $\krylov_m$). Because of the good conjugation properties of basis $(w_i)$ the optimization is decoupled and only one scalar coefficient $\alpha_m$ is sought for at each iteration so that $\|\glam-\glam_{m-1}-\alpha_m w_m\|_\dgschur$ is minimized.

\section{Modification of the FETI method for periodic structures}
For pedagogical reasons the method is first presented in the simplest case, namely a structure made of a single
pattern repeated in a periodic manner; more complex situations are considered in the section following this one.
\subsection{Principles of the method}
\begin{figure}[ht]\centering
\scalebox{.8}{
\ifx\du\undefined
  \newlength{\du}
\fi
\setlength{\du}{15\unitlength}
\begin{tikzpicture}
\pgftransformxscale{1.000000}
\pgftransformyscale{-1.000000}
\definecolor{dialinecolor}{rgb}{0.000000, 0.000000, 0.000000}
\pgfsetstrokecolor{dialinecolor}
\definecolor{dialinecolor}{rgb}{1.000000, 1.000000, 1.000000}
\pgfsetfillcolor{dialinecolor}
\definecolor{dialinecolor}{rgb}{1.000000, 1.000000, 1.000000}
\pgfsetfillcolor{dialinecolor}
\pgfpathellipse{\pgfpoint{13.500000\du}{3.500000\du}}{\pgfpoint{4.500000\du}{0\du}}{\pgfpoint{0\du}{4.500000\du}}
\pgfusepath{fill}
\pgfsetlinewidth{0.100000\du}
\pgfsetdash{}{0pt}
\pgfsetdash{}{0pt}
\definecolor{dialinecolor}{rgb}{0.000000, 0.000000, 0.000000}
\pgfsetstrokecolor{dialinecolor}
\pgfpathellipse{\pgfpoint{13.500000\du}{3.500000\du}}{\pgfpoint{4.500000\du}{0\du}}{\pgfpoint{0\du}{4.500000\du}}
\pgfusepath{stroke}
\definecolor{dialinecolor}{rgb}{0.749020, 0.749020, 0.749020}
\pgfsetfillcolor{dialinecolor}
\pgfpathellipse{\pgfpoint{13.500000\du}{3.500000\du}}{\pgfpoint{2.500000\du}{0\du}}{\pgfpoint{0\du}{2.500000\du}}
\pgfusepath{fill}
\pgfsetlinewidth{0.100000\du}
\pgfsetdash{}{0pt}
\pgfsetdash{}{0pt}
\definecolor{dialinecolor}{rgb}{0.000000, 0.000000, 0.000000}
\pgfsetstrokecolor{dialinecolor}
\pgfpathellipse{\pgfpoint{13.500000\du}{3.500000\du}}{\pgfpoint{2.500000\du}{0\du}}{\pgfpoint{0\du}{2.500000\du}}
\pgfusepath{stroke}
\pgfsetlinewidth{0.100000\du}
\pgfsetdash{}{0pt}
\pgfsetdash{}{0pt}
\pgfsetbuttcap
{
\definecolor{dialinecolor}{rgb}{0.000000, 0.000000, 0.000000}
\pgfsetfillcolor{dialinecolor}
\definecolor{dialinecolor}{rgb}{0.000000, 0.000000, 0.000000}
\pgfsetstrokecolor{dialinecolor}
\draw (13.500000\du,-1.000000\du)--(13.500000\du,1.000000\du);
}
\pgfsetlinewidth{0.100000\du}
\pgfsetdash{}{0pt}
\pgfsetdash{}{0pt}
\pgfsetbuttcap
{
\definecolor{dialinecolor}{rgb}{0.000000, 0.000000, 0.000000}
\pgfsetfillcolor{dialinecolor}
\definecolor{dialinecolor}{rgb}{0.000000, 0.000000, 0.000000}
\pgfsetstrokecolor{dialinecolor}
\draw (9.550000\du,5.687500\du)--(11.250000\du,4.687500\du);
}
\pgfsetlinewidth{0.100000\du}
\pgfsetdash{}{0pt}
\pgfsetdash{}{0pt}
\pgfsetbuttcap
{
\definecolor{dialinecolor}{rgb}{0.000000, 0.000000, 0.000000}
\pgfsetfillcolor{dialinecolor}
\definecolor{dialinecolor}{rgb}{0.000000, 0.000000, 0.000000}
\pgfsetstrokecolor{dialinecolor}
\draw (15.600000\du,4.737500\du)--(17.400000\du,5.687500\du);
}
\definecolor{dialinecolor}{rgb}{0.000000, 0.000000, 0.000000}
\pgfsetstrokecolor{dialinecolor}
\node[anchor=west] at (16.100000\du,-0.212500\du){$\partial^n\Omega$};
\definecolor{dialinecolor}{rgb}{0.000000, 0.000000, 0.000000}
\pgfsetstrokecolor{dialinecolor}
\node[anchor=west] at (9.850000\du,2.287500\du){$\Omega^{(1)}$};
\definecolor{dialinecolor}{rgb}{0.000000, 0.000000, 0.000000}
\pgfsetstrokecolor{dialinecolor}
\node[anchor=west] at (15.950000\du,2.037500\du){$\Omega^{(2)}$};
\definecolor{dialinecolor}{rgb}{0.000000, 0.000000, 0.000000}
\pgfsetstrokecolor{dialinecolor}
\node[anchor=west] at (12.900000\du,7.387500\du){$\Omega^{(3)}$};
\definecolor{dialinecolor}{rgb}{0.000000, 0.000000, 0.000000}
\pgfsetstrokecolor{dialinecolor}
\node[anchor=west] at (11.500000\du,0.087500\du){$\partial^{(1,2)}$};
\definecolor{dialinecolor}{rgb}{0.000000, 0.000000, 0.000000}
\pgfsetstrokecolor{dialinecolor}
\node[anchor=west] at (13.400000\du,0.737500\du){$\partial^{(2,1)}$};
\definecolor{dialinecolor}{rgb}{0.000000, 0.000000, 0.000000}
\pgfsetstrokecolor{dialinecolor}
\node[anchor=west] at (15.800000\du,4.687500\du){$\partial^{(2,3)}$};
\definecolor{dialinecolor}{rgb}{0.000000, 0.000000, 0.000000}
\pgfsetstrokecolor{dialinecolor}
\node[anchor=west] at (15.450000\du,5.887500\du){$\partial^{(3,2)}$};
\definecolor{dialinecolor}{rgb}{0.000000, 0.000000, 0.000000}
\pgfsetstrokecolor{dialinecolor}
\node[anchor=west] at (9.650000\du,5.887500\du){$\partial^{(3,1)}$};
\definecolor{dialinecolor}{rgb}{0.000000, 0.000000, 0.000000}
\pgfsetstrokecolor{dialinecolor}
\node[anchor=west] at (9.100000\du,4.787500\du){$\partial^{(1,3)}$};
\definecolor{dialinecolor}{rgb}{0.000000, 0.000000, 0.000000}
\pgfsetstrokecolor{dialinecolor}
\node[anchor=west] at (13.000000\du,5.487500\du){$\partial^d\Omega$};
\end{tikzpicture}}\caption{``Donut'' toy problem}\label{fig:donut:0}
\end{figure}
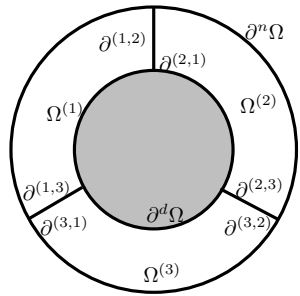

In this toy problem, we consider a periodic structure. In order to introduce concepts without excessive technical
difficulties we first focus on scalar problems so that no local frame is required. All substructures thus have the same
local Schur complement $\schur\s=\schur\n$.
To illustrate our discussion, one can think of the toy problem of determining the temperature field
on a ring whose internal face is submitted to a given temperature (so that no zero energy mode exist in our problem),
the external face being submitted to non-periodic heating.
To fix the ideas, we assume the pattern to be a third of the structure, see figure \ref{fig:donut:0}.

The main assumption underlying the method is that \emph{the numerical information generated at the interface
of one occurrence is pertinent for all identical occurrences in the system}. So, given a search direction $w_i$ and
instead of ``simply'' finding optimal length $\alpha$, a block of search direction $\muv{w}_i$ is generated
by permutation of the information. Referring to figure \ref{fig:donut:0} for the definition of interfaces,
\begin{equation}
w_i =
\begin{pmatrix}
w^{(1,2)} \\
w^{(2,3)} \\
w^{(3,1)}
\end{pmatrix}
\qquad\longrightarrow\qquad
\muv{w}_i =
\begin{pmatrix}
w^{(1,2)} & w^{(2,3)} & w^{(3,1)}\\
w^{(2,3)} & w^{(3,1)} & w^{(1,2)} \\
w^{(3,1)} & w^{(1,2)} & w^{(2,3)}
\end{pmatrix}
\end{equation}
The undertilde notation is used to denote the block of permutations.
As can be seen three different directions of descent have been constructed
by simple permutation of the interface partitions of the preconditioned residual $w_i$.
Then optimization can be realized simultaneously on these three search directions and an optimal linear
combination of these search directions (represented by vector $\uu{\alpha}$) can be deduced so that
$x_{i+1}=x_i +\muv{w}_i\uu{\alpha}_i$. Obviously the search directions $\muv{w}_{i+1}$ for the next iteration
have to be made orthogonal to previous.

Depending on the loading and the structure, the method may converge much faster than classical methods.
In any case the convergence can not be deteriorated because the initial information is still present in the
multivector basis.
The algorithm is presented in the right column of table \ref{tab:alggc}. As can be seen it is quite similar to the
classical conjugate gradient since only a few more operations are involved. The following section discusses
how an efficient implementation can be obtained from the pattern-based framework.

It is to be noted that the algorithm presented in the right column of table \ref{tab:alggc}
applies independently of the way the multivector is generated and thus, as shown in section \ref{sec:extensions},
it can be applied to non-fully-periodic structures. In the case of truly-periodic structures (namely where the interface
trace operators ${\dassem}\s$ are cyclic) computational costs
can be minimized when performing products and when orthogonalizing the new research direction with respect
to previous ones: indeed, because $\dgschur$ and  $\tilde{\dgschur}^{-1}$ themselves exhibit periodicity properties,
the multivectors $\muv{p}_j$ and $\muv{z}_{j+1}$ can also be deduced from one single vector.
For instance, it is sufficient to compute $(\text{line}\ 5.1)\ p_j=Q^T\dgschur w_j$ then to deduce $\muv{p}_j$ and $\muv{w}_j$; it is also  sufficient to compute $(\text{line}\ 5.5)\ z_{j+1}=\tilde{\dgschur}^{-1}r_{j+1}$ then to make $z_{j+1}$ orthogonal to previous multivector of search directions $\muv{w}_{j}$ ($(\text{line}\ 5.6)\ \uu{\beta}_j=-({\muv{w}_j}^T\muv{p}_j)^{-1}(\muv{p}_j^T{z_{j+1}})$ and $(\text{line}\ 5.7)\ w_{j+1}=z_{j+1}+\muv{w}_j\uu{\beta}_j $) and to finally deduce $\muv{w}_{j+1}$ through permutations.

In the following section we discuss how the computations in the multivector FETI algorithm can be organized efficiently
for periodic structures.

\subsection{Efficient implementation of the blocked algorithm}
It is first to be noticed that our method introduces nearly no additional computations compared to the basic method. The main
differences are that
dot-products are replaced by matrix-products $\left({\muv{p}_j}^T\muv{w}_j\right)$ and that the
resulting square matrix needs to be factorized. Moreover, giving the pattern a predominant role
leads to the rewriting of all the operations in a computationally more efficient way: the main idea
is that operations will be realized simultaneously on blocks of vectors (instead of on single vectors).

The principle is that the complete interface can be observed from the point of view of the pattern
with each occurrence of the pattern describing one part of the complete interface. Classically the
interface vector $w$ is stored as subdomain contributions $w^{(s)}$ themselves split in interface contributions
$w^{(s,i)}$. Here however the complete interface is stored as a matrix, named $W$, of local contributions (refer to figure
\ref{fig:donut:0}):
\begin{equation}\label{eq:W}
w=\begin{pmatrix}
w^{(1,2)}\\w^{(2,3)}\\w^{(3,1)}
\end{pmatrix}
\qquad\qquad
w^{(1)}=
\begin{pmatrix}
w^{(1,2)} \\ w^{(1,3)}
\end{pmatrix}\qquad\qquad
W=
\begin{pmatrix}
w^{(1,2)} & w^{(2,3)} & w^{(3,1)} \\ w^{(1,3)} & w^{(2,1)} & w^{(3,2)}
\end{pmatrix}
\end{equation}
where we define $w^{(i,j)}=-w^{(j,i)}$.
Clearly $W$ stores the same information as $w$. The redundancy inside $W$, i.e.\ storing
$w^{(i,j)}$ and $w^{(j,i)}$ corresponds simply to storing the information on different subdomains as is usually done
anyway in parallel implementations when organizing the communications.
In the following capital letters refer to quantities stored in this specific ``block'' format.

Let $\schur\n$ be the Schur complement of the pattern which, in the case of thermal problems,
is identical for all its occurrences. Computing a product $\dgschur w =\sum_s \dassem\s {\schur\n}^+ {\dassem\s}^T w$
can be realized by computing $\left.\schur\n\right.^+ W$ and using adapted assembly operations. The interest of
such an implementation is that the solutions for all the right hand sides (i.e. all columns of $W$) are computed
simultaneously, which enables to minimize the cost of memory transfers which represent a significant contribution
to the overall cost in sparse matrix operations.

Let us now analyse the extra cost associated to the minimization of the residue on a set of research directions.
The issue of applying $\dgschur$ and $\tilde{\dgschur}^{-1}$ on multivectors was discussed above.
One must however note that to put the result of the operation in the multivector format as described for $W$ in
(\ref{eq:W}) a slightly modified implementation of the assembly procedure is required.
The matrix operations such as $\muv{p}^T\muv{w}$ can directly be deduced from the computation of $P^T W$
(using again a modified assembly procedure). In this way this operation is hardly more expensive than
simply computing the vector dot product $p^T w$.
Finally the only significant extra-cost is caused by the factorization of matrix $\muv{p}^T\muv{w}$;
since, in practical problems, the number of occurrences of a pattern is not extremely high
the dimension of this matrix is expected to remain small.
The issue related to the possibility of this minimization matrix to become singular
is discussed in following section.

\subsection{Numerical issues}
As said earlier the proposed algorithm optimizes many computations associated to domain decomposition methods,
and only introduces one additional operation: the factorization of matrix $(\muv{p}^T\muv{w})$ required during
the computation of optimal descent coefficients $\uu{\alpha}_j$. It might occur that for some iteration
this matrix became non-invertible or at least very poorly conditioned. The bad conditioning of matrix $(\muv{p}^T\muv{w})$
can be traced back the non full-column-rankness of the block of search directions $\muv{w}$ which is itself
associated to a redundant information between (at least) two interfaces.
In that case the use of a pseudo inverse is sufficient to avoid breakdowns.

It has been observed that this algorithm, due to the fact that it handles several directions simultaneously,
is more sensitive to roundoff errors than classical conjugate gradient schemes.
For that reason full re-orthogonalization of the directions of descent is often used in order to make the iteration
procedure more robust. The $\dgschur$-normalization of search directions is then a simplifying feature;
after the computation of $\uuu{\eta}=(\muv{p}^T\muv{w})$, one substitutes $\muv{w}\leftarrow\muv{w}\uuu{\eta}^{-\frac{1}{2}}$
and $\muv{p}\leftarrow\muv{p}\uuu{\eta}^{-\frac{1}{2}}$ such that $(\muv{p}^T\muv{w})=\id$. As a result no system of equation
needs to be solved during the full-reorthogonalization and residue minimization steps.
Obviously the normalization procedure described above is the right place to handle redundancy of directions of descent in
the multivector.

\subsection{Assessments}
\subsubsection{Academical tests}
\label{subsubsec_aca_per}
\begin{figure}[ht]\centering
   \begin{minipage}[c]{.49\linewidth}\centering
      \includegraphics[width=0.99\textwidth]{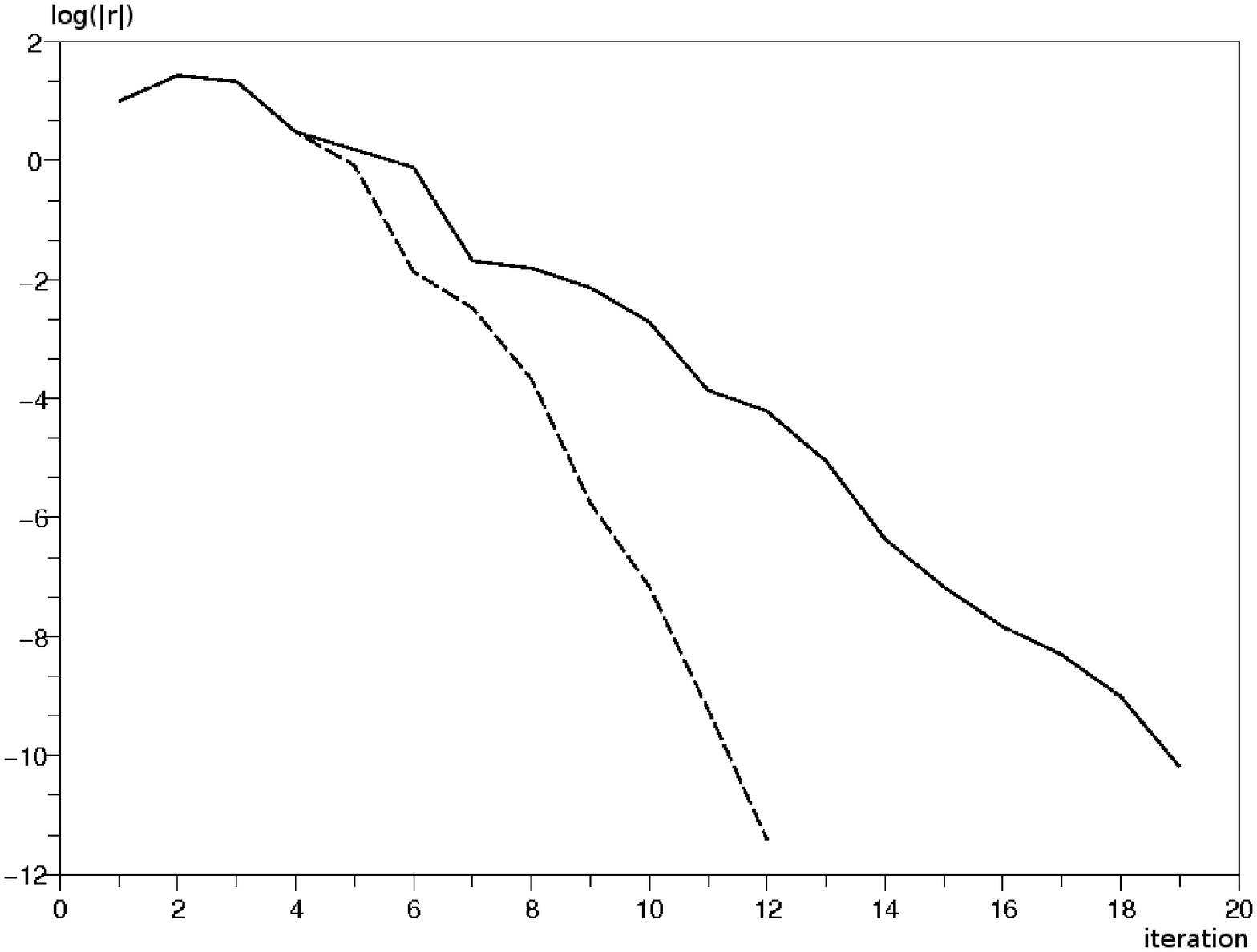}\\Five repetitions
   \end{minipage} \hfill
   \begin{minipage}[c]{.49\linewidth}\centering
      \includegraphics[width=0.99\textwidth]{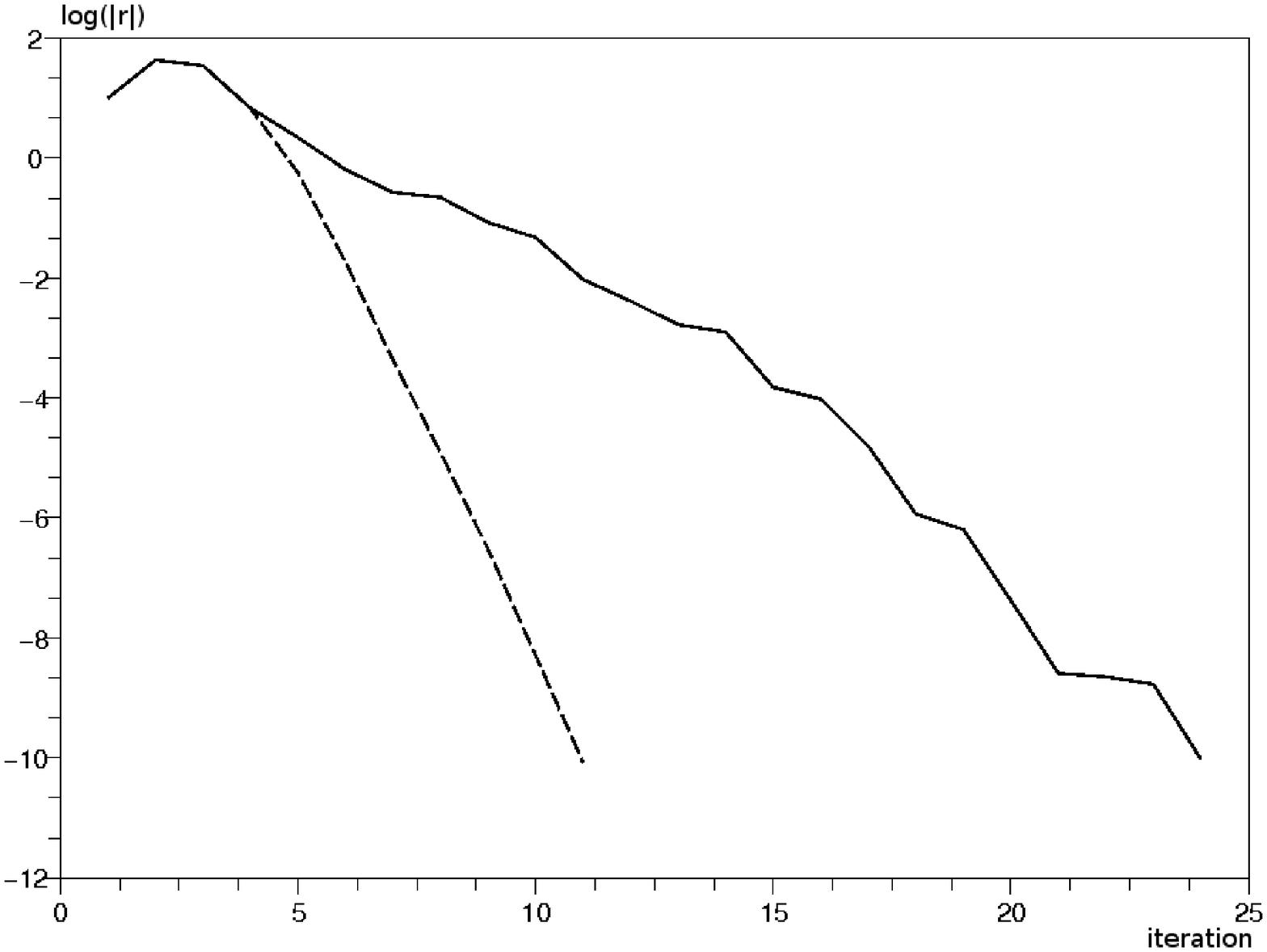}\\Nine repetitions
   \end{minipage}\caption{Comparison between new (dashed) and classical (solid) algorithms}\label{fig:aca_result:0}
\end{figure}

Here the Schur complement of the pattern is represented by a random dense matrix made symmetric definite positive.
The Dirichlet's preconditioner is used. The pattern is repeated either five or nine times. The loading is random and non-periodic.
Figure \ref{fig:aca_result:0} presents the number of iterations of the conjugate gradient to reach a given precision
for a classical dual domain decomposition method and for the block (multivector) version (implemented in Scilab,
\url{www.scilab.org}).
As can be seen, using several directions of descent per iteration can reduce the number of iteration
by about $35\%$  (respectively $55\%$) for the five-repetition example (respectively for the nine-repetition test), reminding
that the cost per iteration is nearly identical.

As can be observed in almost every cases, the norms of the residual during the first iterations of the block and non-block methods are almost coincident. 
Indeed, the first iterations of conjugate gradient explore the higher part of the active spectrum of the problem which corresponds to large wavelength phenomena which are almost identical in any column of the multivector. Once these contributions are found, the enrichment of the search space offered by multivectors becomes effective and leads to significant improvement of the convergence.


\subsubsection{Donut tests}
\begin{figure}[ht]\centering
   \begin{minipage}[c]{.20\linewidth}\centering
      \includegraphics[width=0.98\textwidth]{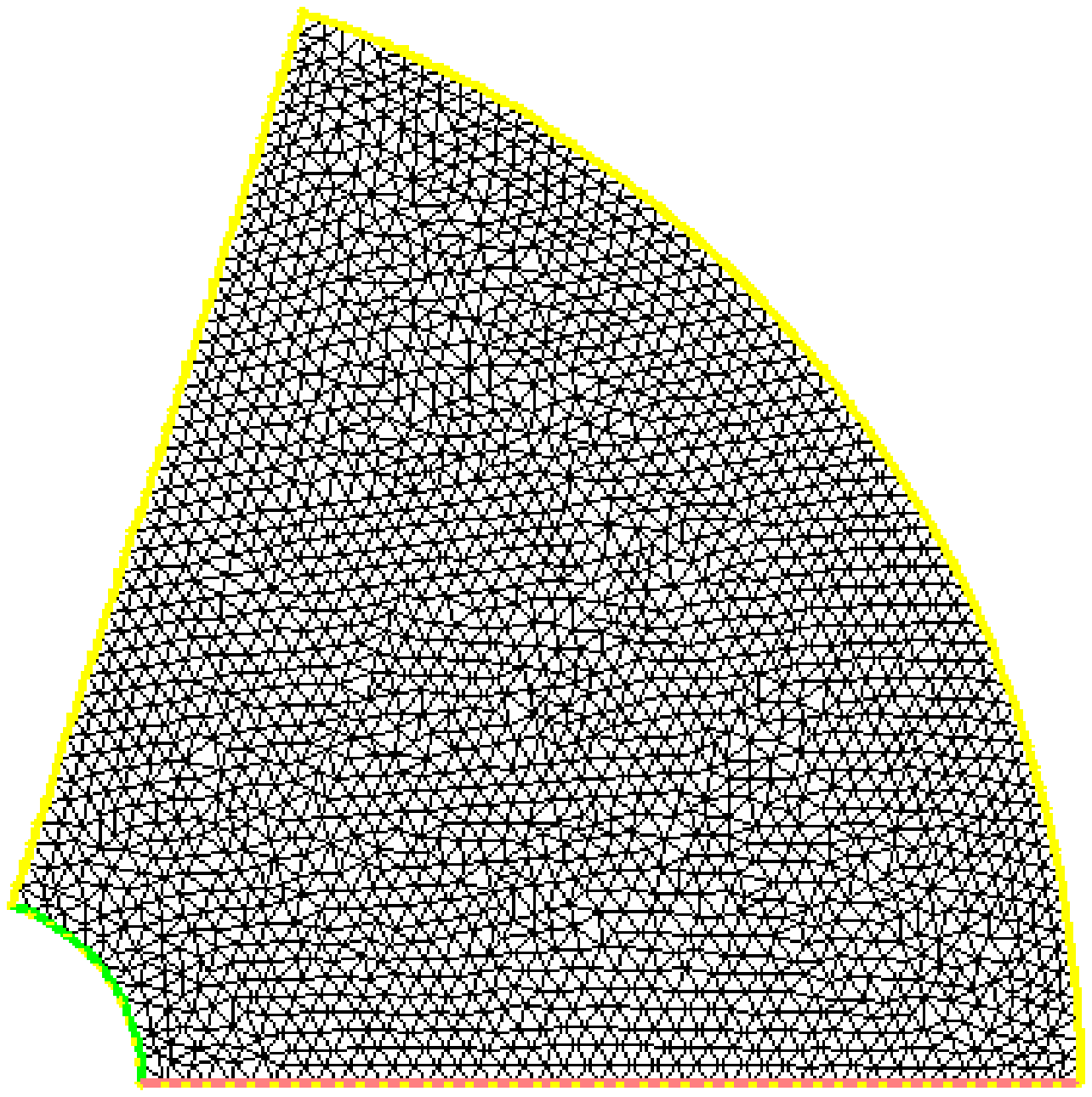}\\Pattern \mbox{(Five-part donut)}
   \end{minipage} \hfill
   \begin{minipage}[c]{.39\linewidth}\centering
      \includegraphics[width=0.99\textwidth]{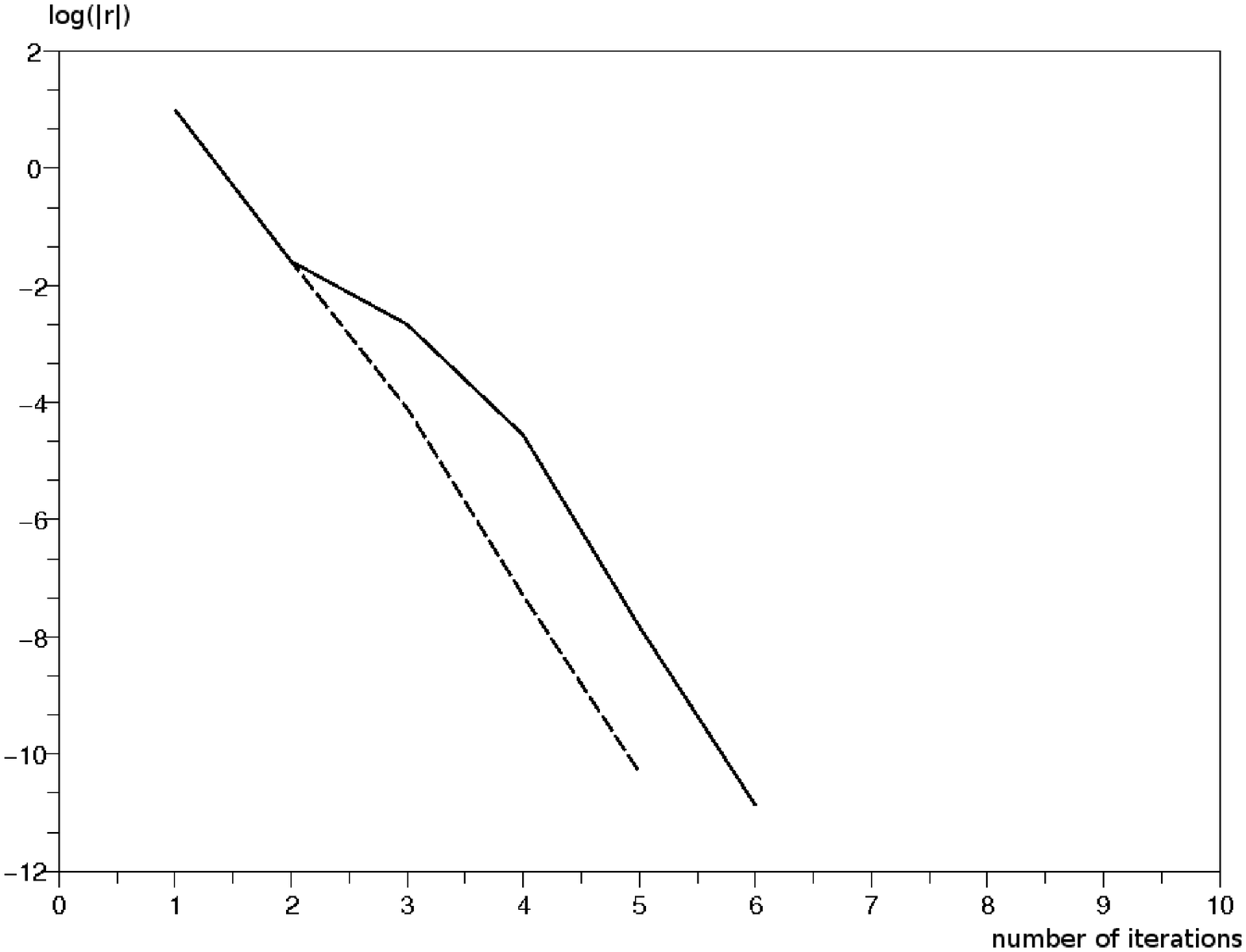}\\Five-part donut
   \end{minipage} \hfill
   \begin{minipage}[c]{.39\linewidth}\centering
      \includegraphics[width=0.99\textwidth]{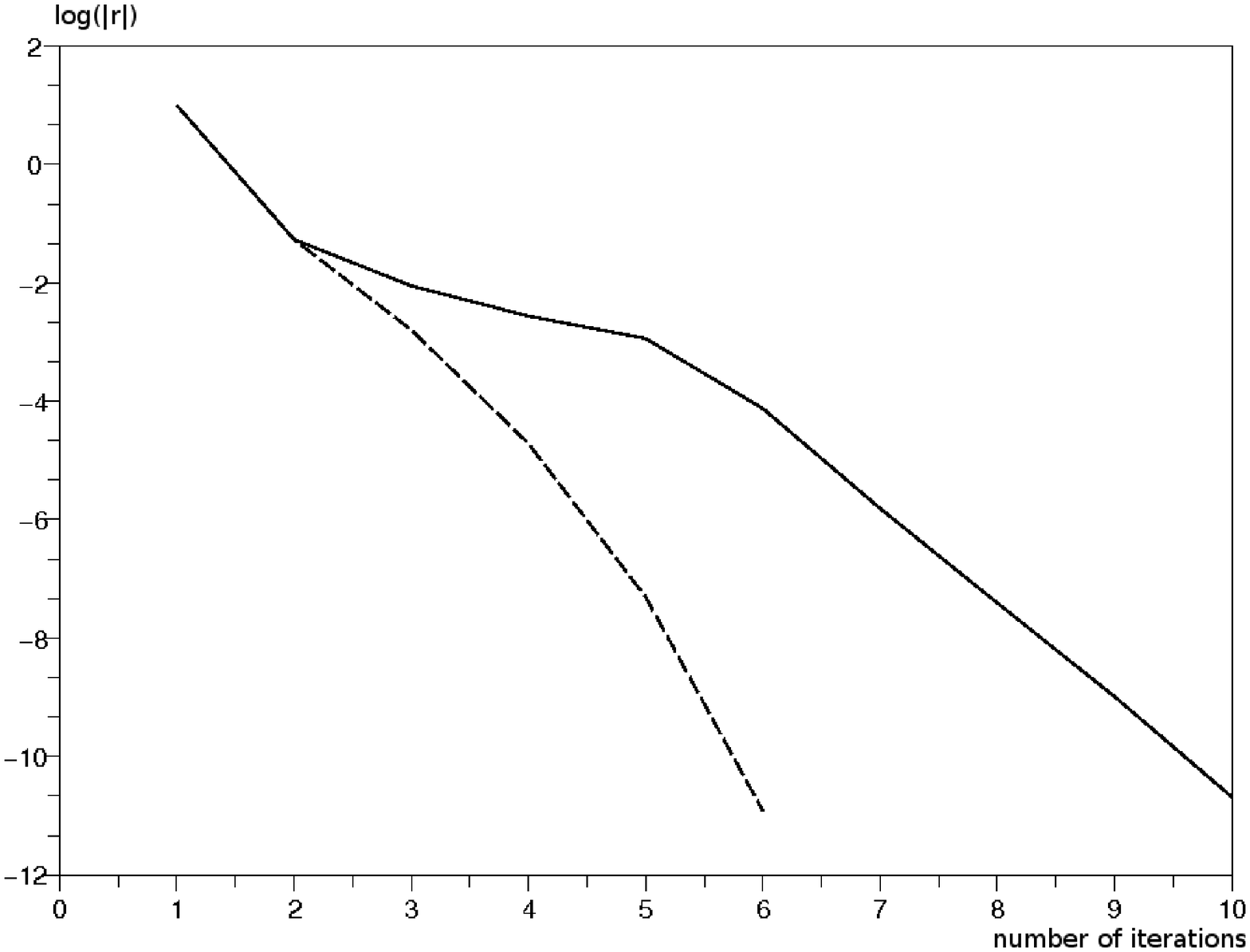}\\Nine-part donut
   \end{minipage}\caption{Comparison between new (dashed) and classical (solid) algorithms}\label{fig:donut_result:0}
\end{figure}

In order to evaluate the performance results on more realistic cases the conductivity matrix of a piece of
donut is obtained from FreeFem++ (\url{www.freefem.org}) and imported in Scilab. The pattern in these examples
has about 2 000 degrees of freedom (dof) of which about 100 are on its interface. Here again the
loading is random and non-periodic. The Dirichlet preconditioner is used. Figure \ref{fig:donut_result:0} presents the number of iterations of
the conjugate gradient to reach a given precision for a classical dual domain decomposition method and for the
block dual domain decomposition method. The algorithms are assessed on the five and nine-part donuts.
The decrease of the number of iterations is less spectacular than in previous cases,
although with 9 repetitions the block algorithm converges in about $50\%$ less iterations.

It is also interesting to estimate the computational costs of the classical FETI solver and the block variant for
systems with repeated patterns (and applied to periodic structure in this section). Here the algorithm is programmed in
scilab (an interpreted programming language) and cpu time are only indicative. Nevertheless it is informative to
compare the CPU time for three different variants of the FETI solver:
\begin{itemize}
\item A basic (``classical'') FETI method where the fact that the structure is made of repeated patterns is not at all
accounted for.
\item The same classical FETI algorithm but where the periodicity is utilized to compute efficiently
$\dgschur w =\sum_s \dassem\s {\schur\n}^+ {\dassem\s}^T w$. Since every subdomain has the same local operator all
subdomain contributions in this interface problem are computed simultaneously, namely applying the forward and backward substitution
at once for several right hand sides. Clearly this algorithm will exhibit the same convergence history but the cost related to the factorization of $\schur\n$ is paid only once, and the cost per iteration is significantly reduced. This version of the FETI method will be called the classical FETI-mrhs.
\item
Finally we apply the multivector FETI algorithm as proposed in this paper where several directions of descent are considered
per iteration whereas the cost per iteration is similar to the FETI method denoted above as the classical FETI-mrhs.
\end{itemize}

The cpu times are listed in Table \ref{tab:donut1}.
Obviously the performance of  the classical FETI is significantly penalized by the fact that no use is made of the
property that local operators are identical; nevertheless if implemented on a multiprocessor computer and
assuming perfect parallelism (i.e.\ dividing the cpu time for the classical FETI by the number of domains) the CPU time
would be very close to the one of the classical FETI-mrhs. Hence in a parallel computing setting using the fact that
local operators are identical as in the classical FETI-mrhs approach is not very advantageous.
When applying the multivector FETI however faster computing will be achieved thanks to the fact
that the entire information generated at every iteration is exploited in the multiple search directions,
thereby leading to a faster convergence with a similar cost per iteration. It amounts to a gain
of 50\% in computing time in our test for the 9-part donut case.
\begin{table}[ht]\centering
 \begin{tabular}{|c|p{3cm}|p{3.1cm}|p{3cm}|}\hline
\multicolumn{4}{|c|}{5-part donut}\\\hline
& classical FETI & classical FETI-mrhs & multivector FETI   \\\hline
number of iterations &5&5&4\\  \hline
CPU time (s) &31.04&6.41&6.37\\\hline\hline
\multicolumn{4}{|c|}{9-part donut}\\\hline
& classical FETI & classical FETI-mrhs & multivector FETI \\\hline
number of iterations &9&9&5\\  \hline
CPU time (s) &72.6&8.49&4.39\\\hline
 \end{tabular}\caption{Performance of the multivector approach for thermal donut problems}\label{tab:donut1}
\end{table}

\section{Extensions}\label{sec:extensions}
\subsection{Non-fully periodic structures}
\begin{figure}[ht]\centering
\subfigure[Donut with one stand]{\scalebox{.75}{
\ifx\du\undefined
  \newlength{\du}
\fi
\setlength{\du}{15\unitlength}
\begin{tikzpicture}
\pgftransformxscale{1.000000}
\pgftransformyscale{-1.000000}
\definecolor{dialinecolor}{rgb}{0.000000, 0.000000, 0.000000}
\pgfsetstrokecolor{dialinecolor}
\definecolor{dialinecolor}{rgb}{1.000000, 1.000000, 1.000000}
\pgfsetfillcolor{dialinecolor}
\pgfsetlinewidth{0.100000\du}
\pgfsetdash{}{0pt}
\pgfsetdash{}{0pt}
\pgfsetmiterjoin
\definecolor{dialinecolor}{rgb}{1.000000, 1.000000, 1.000000}
\pgfsetfillcolor{dialinecolor}
\fill (10.950000\du,7.087500\du)--(10.950000\du,9.850000\du)--(16.000000\du,9.850000\du)--(16.000000\du,7.087500\du)--cycle;
\definecolor{dialinecolor}{rgb}{0.000000, 0.000000, 0.000000}
\pgfsetstrokecolor{dialinecolor}
\draw (10.950000\du,7.087500\du)--(10.950000\du,9.850000\du)--(16.000000\du,9.850000\du)--(16.000000\du,7.087500\du)--cycle;
\definecolor{dialinecolor}{rgb}{1.000000, 1.000000, 1.000000}
\pgfsetfillcolor{dialinecolor}
\pgfpathellipse{\pgfpoint{13.500000\du}{3.500000\du}}{\pgfpoint{4.500000\du}{0\du}}{\pgfpoint{0\du}{4.500000\du}}
\pgfusepath{fill}
\pgfsetlinewidth{0.100000\du}
\pgfsetdash{}{0pt}
\pgfsetdash{}{0pt}
\definecolor{dialinecolor}{rgb}{0.000000, 0.000000, 0.000000}
\pgfsetstrokecolor{dialinecolor}
\pgfpathellipse{\pgfpoint{13.500000\du}{3.500000\du}}{\pgfpoint{4.500000\du}{0\du}}{\pgfpoint{0\du}{4.500000\du}}
\pgfusepath{stroke}
\definecolor{dialinecolor}{rgb}{0.749020, 0.749020, 0.749020}
\pgfsetfillcolor{dialinecolor}
\pgfpathellipse{\pgfpoint{13.500000\du}{3.500000\du}}{\pgfpoint{2.500000\du}{0\du}}{\pgfpoint{0\du}{2.500000\du}}
\pgfusepath{fill}
\pgfsetlinewidth{0.100000\du}
\pgfsetdash{}{0pt}
\pgfsetdash{}{0pt}
\definecolor{dialinecolor}{rgb}{0.000000, 0.000000, 0.000000}
\pgfsetstrokecolor{dialinecolor}
\pgfpathellipse{\pgfpoint{13.500000\du}{3.500000\du}}{\pgfpoint{2.500000\du}{0\du}}{\pgfpoint{0\du}{2.500000\du}}
\pgfusepath{stroke}
\pgfsetlinewidth{0.100000\du}
\pgfsetdash{}{0pt}
\pgfsetdash{}{0pt}
\pgfsetbuttcap
{
\definecolor{dialinecolor}{rgb}{0.000000, 0.000000, 0.000000}
\pgfsetfillcolor{dialinecolor}
\definecolor{dialinecolor}{rgb}{0.000000, 0.000000, 0.000000}
\pgfsetstrokecolor{dialinecolor}
\draw (13.500000\du,-1.000000\du)--(13.500000\du,1.000000\du);
}
\pgfsetlinewidth{0.100000\du}
\pgfsetdash{}{0pt}
\pgfsetdash{}{0pt}
\pgfsetbuttcap
{
\definecolor{dialinecolor}{rgb}{0.000000, 0.000000, 0.000000}
\pgfsetfillcolor{dialinecolor}
\definecolor{dialinecolor}{rgb}{0.000000, 0.000000, 0.000000}
\pgfsetstrokecolor{dialinecolor}
\draw (9.550000\du,5.687500\du)--(11.250000\du,4.687500\du);
}
\pgfsetlinewidth{0.100000\du}
\pgfsetdash{}{0pt}
\pgfsetdash{}{0pt}
\pgfsetbuttcap
{
\definecolor{dialinecolor}{rgb}{0.000000, 0.000000, 0.000000}
\pgfsetfillcolor{dialinecolor}
\definecolor{dialinecolor}{rgb}{0.000000, 0.000000, 0.000000}
\pgfsetstrokecolor{dialinecolor}
\draw (15.600000\du,4.737500\du)--(17.400000\du,5.687500\du);
}
\definecolor{dialinecolor}{rgb}{0.000000, 0.000000, 0.000000}
\pgfsetstrokecolor{dialinecolor}
\node[anchor=west] at (9.650000\du,2.650000\du){$\Omega^{(1)}$};
\definecolor{dialinecolor}{rgb}{0.000000, 0.000000, 0.000000}
\pgfsetstrokecolor{dialinecolor}
\node[anchor=west] at (16.250000\du,2.600000\du){$\Omega^{(2)}$};
\definecolor{dialinecolor}{rgb}{0.000000, 0.000000, 0.000000}
\pgfsetstrokecolor{dialinecolor}
\node[anchor=west] at (13.050000\du,6.950000\du){$\Omega^{(3)}$};
\definecolor{dialinecolor}{rgb}{0.000000, 0.000000, 0.000000}
\pgfsetstrokecolor{dialinecolor}
\node[anchor=west] at (12.650000\du,9.150000\du){$\Omega^{(4)}$};
\end{tikzpicture}}}\hfill
\subfigure[Donut with two stands]{\scalebox{.75}{
\ifx\du\undefined
  \newlength{\du}
\fi
\setlength{\du}{15\unitlength}
\begin{tikzpicture}
\pgftransformxscale{1.000000}
\pgftransformyscale{-1.000000}
\definecolor{dialinecolor}{rgb}{0.000000, 0.000000, 0.000000}
\pgfsetstrokecolor{dialinecolor}
\definecolor{dialinecolor}{rgb}{1.000000, 1.000000, 1.000000}
\pgfsetfillcolor{dialinecolor}
\pgfsetlinewidth{0.100000\du}
\pgfsetdash{}{0pt}
\pgfsetdash{}{0pt}
\pgfsetmiterjoin
\definecolor{dialinecolor}{rgb}{1.000000, 1.000000, 1.000000}
\pgfsetfillcolor{dialinecolor}
\fill (7.000000\du,0.000000\du)--(7.000000\du,2.762500\du)--(12.050000\du,2.762500\du)--(12.050000\du,0.000000\du)--cycle;
\definecolor{dialinecolor}{rgb}{0.000000, 0.000000, 0.000000}
\pgfsetstrokecolor{dialinecolor}
\draw (7.000000\du,0.000000\du)--(7.000000\du,2.762500\du)--(12.050000\du,2.762500\du)--(12.050000\du,0.000000\du)--cycle;
\pgfsetlinewidth{0.100000\du}
\pgfsetdash{}{0pt}
\pgfsetdash{}{0pt}
\pgfsetmiterjoin
\definecolor{dialinecolor}{rgb}{1.000000, 1.000000, 1.000000}
\pgfsetfillcolor{dialinecolor}
\fill (15.000000\du,0.000000\du)--(15.000000\du,2.762500\du)--(20.050000\du,2.762500\du)--(20.050000\du,0.000000\du)--cycle;
\definecolor{dialinecolor}{rgb}{0.000000, 0.000000, 0.000000}
\pgfsetstrokecolor{dialinecolor}
\draw (15.000000\du,0.000000\du)--(15.000000\du,2.762500\du)--(20.050000\du,2.762500\du)--(20.050000\du,0.000000\du)--cycle;
\definecolor{dialinecolor}{rgb}{1.000000, 1.000000, 1.000000}
\pgfsetfillcolor{dialinecolor}
\pgfpathellipse{\pgfpoint{13.500000\du}{3.500000\du}}{\pgfpoint{4.500000\du}{0\du}}{\pgfpoint{0\du}{4.500000\du}}
\pgfusepath{fill}
\pgfsetlinewidth{0.100000\du}
\pgfsetdash{}{0pt}
\pgfsetdash{}{0pt}
\definecolor{dialinecolor}{rgb}{0.000000, 0.000000, 0.000000}
\pgfsetstrokecolor{dialinecolor}
\pgfpathellipse{\pgfpoint{13.500000\du}{3.500000\du}}{\pgfpoint{4.500000\du}{0\du}}{\pgfpoint{0\du}{4.500000\du}}
\pgfusepath{stroke}
\definecolor{dialinecolor}{rgb}{0.749020, 0.749020, 0.749020}
\pgfsetfillcolor{dialinecolor}
\pgfpathellipse{\pgfpoint{13.500000\du}{3.500000\du}}{\pgfpoint{2.500000\du}{0\du}}{\pgfpoint{0\du}{2.500000\du}}
\pgfusepath{fill}
\pgfsetlinewidth{0.100000\du}
\pgfsetdash{}{0pt}
\pgfsetdash{}{0pt}
\definecolor{dialinecolor}{rgb}{0.000000, 0.000000, 0.000000}
\pgfsetstrokecolor{dialinecolor}
\pgfpathellipse{\pgfpoint{13.500000\du}{3.500000\du}}{\pgfpoint{2.500000\du}{0\du}}{\pgfpoint{0\du}{2.500000\du}}
\pgfusepath{stroke}
\pgfsetlinewidth{0.100000\du}
\pgfsetdash{}{0pt}
\pgfsetdash{}{0pt}
\pgfsetbuttcap
{
\definecolor{dialinecolor}{rgb}{0.000000, 0.000000, 0.000000}
\pgfsetfillcolor{dialinecolor}
\definecolor{dialinecolor}{rgb}{0.000000, 0.000000, 0.000000}
\pgfsetstrokecolor{dialinecolor}
\draw (13.500000\du,-1.000000\du)--(13.500000\du,1.000000\du);
}
\pgfsetlinewidth{0.100000\du}
\pgfsetdash{}{0pt}
\pgfsetdash{}{0pt}
\pgfsetbuttcap
{
\definecolor{dialinecolor}{rgb}{0.000000, 0.000000, 0.000000}
\pgfsetfillcolor{dialinecolor}
\definecolor{dialinecolor}{rgb}{0.000000, 0.000000, 0.000000}
\pgfsetstrokecolor{dialinecolor}
\draw (9.550000\du,5.687500\du)--(11.250000\du,4.687500\du);
}
\pgfsetlinewidth{0.100000\du}
\pgfsetdash{}{0pt}
\pgfsetdash{}{0pt}
\pgfsetbuttcap
{
\definecolor{dialinecolor}{rgb}{0.000000, 0.000000, 0.000000}
\pgfsetfillcolor{dialinecolor}
\definecolor{dialinecolor}{rgb}{0.000000, 0.000000, 0.000000}
\pgfsetstrokecolor{dialinecolor}
\draw (15.600000\du,4.737500\du)--(17.400000\du,5.687500\du);
}
\definecolor{dialinecolor}{rgb}{0.000000, 0.000000, 0.000000}
\pgfsetstrokecolor{dialinecolor}
\node[anchor=west] at (9.650000\du,2.650000\du){$\Omega^{(1)}$};
\definecolor{dialinecolor}{rgb}{0.000000, 0.000000, 0.000000}
\pgfsetstrokecolor{dialinecolor}
\node[anchor=west] at (16.250000\du,2.600000\du){$\Omega^{(2)}$};
\definecolor{dialinecolor}{rgb}{0.000000, 0.000000, 0.000000}
\pgfsetstrokecolor{dialinecolor}
\node[anchor=west] at (13.050000\du,6.950000\du){$\Omega^{(3)}$};
\definecolor{dialinecolor}{rgb}{0.000000, 0.000000, 0.000000}
\pgfsetstrokecolor{dialinecolor}
\node[anchor=west] at (18.300000\du,1.550000\du){$\Omega^{(4)}$};
\definecolor{dialinecolor}{rgb}{0.000000, 0.000000, 0.000000}
\pgfsetstrokecolor{dialinecolor}
\node[anchor=west] at (7.600000\du,1.550000\du){$\Omega^{(5)}$};
\end{tikzpicture}}}\hfill
\subfigure[Non periodic Dirichlet bc's]{\scalebox{.75}{
\ifx\du\undefined
  \newlength{\du}
\fi
\setlength{\du}{15\unitlength}
\begin{tikzpicture}
\pgftransformxscale{1.000000}
\pgftransformyscale{-1.000000}
\definecolor{dialinecolor}{rgb}{0.000000, 0.000000, 0.000000}
\pgfsetstrokecolor{dialinecolor}
\definecolor{dialinecolor}{rgb}{1.000000, 1.000000, 1.000000}
\pgfsetfillcolor{dialinecolor}
\pgfsetlinewidth{0.100000\du}
\pgfsetdash{}{0pt}
\pgfsetdash{}{0pt}
\pgfsetmiterjoin
{\pgfsetcornersarced{\pgfpoint{0.000000\du}{0.000000\du}}\definecolor{dialinecolor}{rgb}{1.000000, 1.000000, 1.000000}
\pgfsetfillcolor{dialinecolor}
\fill (11.095000\du,7.035000\du)--(11.095000\du,9.797500\du)--(16.145000\du,9.797500\du)--(16.145000\du,7.035000\du)--cycle;
}{\pgfsetcornersarced{\pgfpoint{0.000000\du}{0.000000\du}}\definecolor{dialinecolor}{rgb}{1.000000, 1.000000, 1.000000}
\pgfsetstrokecolor{dialinecolor}
\draw (11.095000\du,7.035000\du)--(11.095000\du,9.797500\du)--(16.145000\du,9.797500\du)--(16.145000\du,7.035000\du)--cycle;
}\definecolor{dialinecolor}{rgb}{1.000000, 1.000000, 1.000000}
\pgfsetfillcolor{dialinecolor}
\pgfpathellipse{\pgfpoint{13.450000\du}{4.000000\du}}{\pgfpoint{4.600000\du}{0\du}}{\pgfpoint{0\du}{4.200000\du}}
\pgfusepath{fill}
\pgfsetlinewidth{0.100000\du}
\pgfsetdash{}{0pt}
\pgfsetdash{}{0pt}
\definecolor{dialinecolor}{rgb}{0.000000, 0.000000, 0.000000}
\pgfsetstrokecolor{dialinecolor}
\pgfpathellipse{\pgfpoint{13.450000\du}{4.000000\du}}{\pgfpoint{4.600000\du}{0\du}}{\pgfpoint{0\du}{4.200000\du}}
\pgfusepath{stroke}
\pgfsetlinewidth{0.100000\du}
\pgfsetdash{}{0pt}
\pgfsetdash{}{0pt}
\pgfsetmiterjoin
{\pgfsetcornersarced{\pgfpoint{0.000000\du}{0.000000\du}}\definecolor{dialinecolor}{rgb}{1.000000, 1.000000, 1.000000}
\pgfsetfillcolor{dialinecolor}
\fill (8.550000\du,0.500000\du)--(8.550000\du,7.150000\du)--(18.350000\du,7.150000\du)--(18.350000\du,0.500000\du)--cycle;
}{\pgfsetcornersarced{\pgfpoint{0.000000\du}{0.000000\du}}\definecolor{dialinecolor}{rgb}{1.000000, 1.000000, 1.000000}
\pgfsetstrokecolor{dialinecolor}
\draw (8.550000\du,0.500000\du)--(8.550000\du,7.150000\du)--(18.350000\du,7.150000\du)--(18.350000\du,0.500000\du)--cycle;
}\definecolor{dialinecolor}{rgb}{1.000000, 1.000000, 1.000000}
\pgfsetfillcolor{dialinecolor}
\pgfpathellipse{\pgfpoint{13.500000\du}{3.500000\du}}{\pgfpoint{4.500000\du}{0\du}}{\pgfpoint{0\du}{4.500000\du}}
\pgfusepath{fill}
\pgfsetlinewidth{0.100000\du}
\pgfsetdash{}{0pt}
\pgfsetdash{}{0pt}
\definecolor{dialinecolor}{rgb}{0.000000, 0.000000, 0.000000}
\pgfsetstrokecolor{dialinecolor}
\pgfpathellipse{\pgfpoint{13.500000\du}{3.500000\du}}{\pgfpoint{4.500000\du}{0\du}}{\pgfpoint{0\du}{4.500000\du}}
\pgfusepath{stroke}
\definecolor{dialinecolor}{rgb}{0.749020, 0.749020, 0.749020}
\pgfsetfillcolor{dialinecolor}
\pgfpathellipse{\pgfpoint{13.500000\du}{3.500000\du}}{\pgfpoint{2.500000\du}{0\du}}{\pgfpoint{0\du}{2.500000\du}}
\pgfusepath{fill}
\pgfsetlinewidth{0.100000\du}
\pgfsetdash{}{0pt}
\pgfsetdash{}{0pt}
\definecolor{dialinecolor}{rgb}{0.000000, 0.000000, 0.000000}
\pgfsetstrokecolor{dialinecolor}
\pgfpathellipse{\pgfpoint{13.500000\du}{3.500000\du}}{\pgfpoint{2.500000\du}{0\du}}{\pgfpoint{0\du}{2.500000\du}}
\pgfusepath{stroke}
\pgfsetlinewidth{0.100000\du}
\pgfsetdash{}{0pt}
\pgfsetdash{}{0pt}
\pgfsetbuttcap
{
\definecolor{dialinecolor}{rgb}{0.000000, 0.000000, 0.000000}
\pgfsetfillcolor{dialinecolor}
\definecolor{dialinecolor}{rgb}{0.000000, 0.000000, 0.000000}
\pgfsetstrokecolor{dialinecolor}
\draw (13.500000\du,-1.000000\du)--(13.500000\du,1.000000\du);
}
\pgfsetlinewidth{0.100000\du}
\pgfsetdash{}{0pt}
\pgfsetdash{}{0pt}
\pgfsetbuttcap
{
\definecolor{dialinecolor}{rgb}{0.000000, 0.000000, 0.000000}
\pgfsetfillcolor{dialinecolor}
\definecolor{dialinecolor}{rgb}{0.000000, 0.000000, 0.000000}
\pgfsetstrokecolor{dialinecolor}
\draw (9.550000\du,5.687500\du)--(11.250000\du,4.687500\du);
}
\pgfsetlinewidth{0.100000\du}
\pgfsetdash{}{0pt}
\pgfsetdash{}{0pt}
\pgfsetbuttcap
{
\definecolor{dialinecolor}{rgb}{0.000000, 0.000000, 0.000000}
\pgfsetfillcolor{dialinecolor}
\definecolor{dialinecolor}{rgb}{0.000000, 0.000000, 0.000000}
\pgfsetstrokecolor{dialinecolor}
\draw (15.600000\du,4.737500\du)--(17.400000\du,5.687500\du);
}
\definecolor{dialinecolor}{rgb}{0.000000, 0.000000, 0.000000}
\pgfsetstrokecolor{dialinecolor}
\node[anchor=west] at (9.400000\du,3.000000\du){$\Omega^{(1)}$};
\definecolor{dialinecolor}{rgb}{0.000000, 0.000000, 0.000000}
\pgfsetstrokecolor{dialinecolor}
\node[anchor=west] at (16.150000\du,2.600000\du){$\Omega^{(2)}$};
\definecolor{dialinecolor}{rgb}{0.000000, 0.000000, 0.000000}
\pgfsetstrokecolor{dialinecolor}
\node[anchor=west] at (13.050000\du,6.950000\du){$\Omega^{(3)}$};
\definecolor{dialinecolor}{rgb}{0.000000, 0.000000, 0.000000}
\pgfsetstrokecolor{dialinecolor}
\node[anchor=west] at (16.800000\du,0.150000\du){$\partial^n\Omega$};
\definecolor{dialinecolor}{rgb}{0.000000, 0.000000, 0.000000}
\pgfsetstrokecolor{dialinecolor}
\node[anchor=west] at (13.000000\du,9.000000\du){$\partial^d\Omega$};
\pgfsetlinewidth{0.100000\du}
\pgfsetdash{}{0pt}
\pgfsetdash{}{0pt}
\pgfsetbuttcap
{
\definecolor{dialinecolor}{rgb}{0.000000, 0.000000, 0.000000}
\pgfsetfillcolor{dialinecolor}
\definecolor{dialinecolor}{rgb}{0.000000, 0.000000, 0.000000}
\pgfsetstrokecolor{dialinecolor}
\draw (13.450000\du,3.825000\du)--(13.450000\du,3.825000\du);
}
\pgfsetlinewidth{0.100000\du}
\pgfsetdash{}{0pt}
\pgfsetdash{}{0pt}
\pgfsetbuttcap
{
\definecolor{dialinecolor}{rgb}{0.000000, 0.000000, 0.000000}
\pgfsetfillcolor{dialinecolor}
\definecolor{dialinecolor}{rgb}{0.000000, 0.000000, 0.000000}
\pgfsetstrokecolor{dialinecolor}
\draw (16.300000\du,7.050000\du)--(16.600000\du,7.250000\du);
}
\pgfsetlinewidth{0.100000\du}
\pgfsetdash{}{0pt}
\pgfsetdash{}{0pt}
\pgfsetbuttcap
{
\definecolor{dialinecolor}{rgb}{0.000000, 0.000000, 0.000000}
\pgfsetfillcolor{dialinecolor}
\definecolor{dialinecolor}{rgb}{0.000000, 0.000000, 0.000000}
\pgfsetstrokecolor{dialinecolor}
\draw (10.400000\du,7.350000\du)--(10.646400\du,7.044050\du);
}
\end{tikzpicture}}}
\caption{Various extensions}\label{fig:donut:n}
\end{figure}
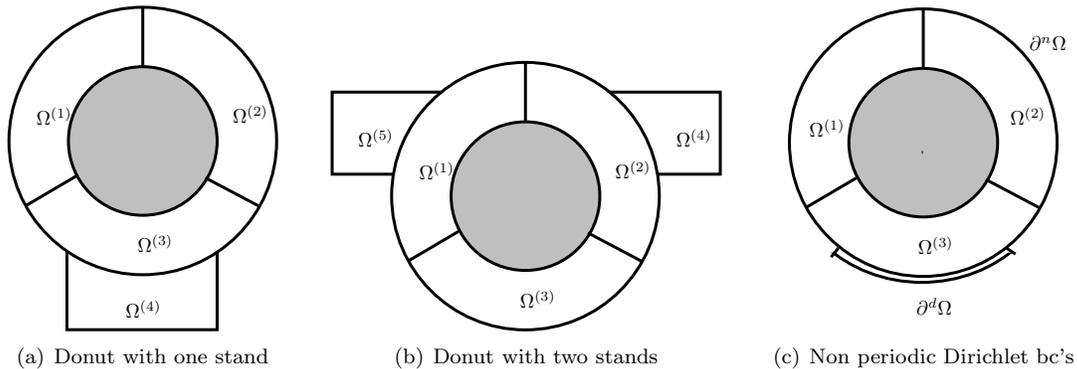

When a structure can not be decomposed into several occurrences of a single pattern, it cannot be treated with the procedure
described in the previous session dealing with periodic structures. It is to be noticed that boundary conditions and/or
connection conditions between
subdomains can also alter the periodicity.
Figure \ref{fig:donut:n} presents different cases where
the strategy developed in the previous section needs to be slightly modified to still be able to make efficient use of
the presence of patterns.
\subsubsection{Donut with one stand}
Let us consider a donut with one stand (subdomain 4) connected to subdomain 3 (figure \ref{fig:donut:n}.a).
There are thus two patterns, one of which being instantiated three times. The topology of interfaces is modified and
the system cannot be seen as periodic anymore. Nevertheless one can still build a multivector that makes
use of the information computed on the interfaces between repeating subdomains. The simplest multivector one
can consider is:
\begin{equation}
w =
\begin{pmatrix}
w^{(1,2)} \\
w^{(2,3)} \\
w^{(3,1)} \\
w^{(3,4)}
\end{pmatrix}
\qquad
\muv{w} =
\begin{pmatrix}
w^{(1,2)} & w^{(2,3)} & w^{(3,1)} \\
w^{(2,3)} & w^{(3,1)} & w^{(1,2)} \\
w^{(3,1)} & w^{(1,2)} & w^{(2,3)} \\
w^{(3,4)} & w^{(3,4)} & w^{(3,4)}
\end{pmatrix}
\end{equation}

\begin{figure}[ht]\centering
   \begin{minipage}[c]{.47\linewidth}\centering
      \includegraphics[angle=-90,width=0.99\textwidth]{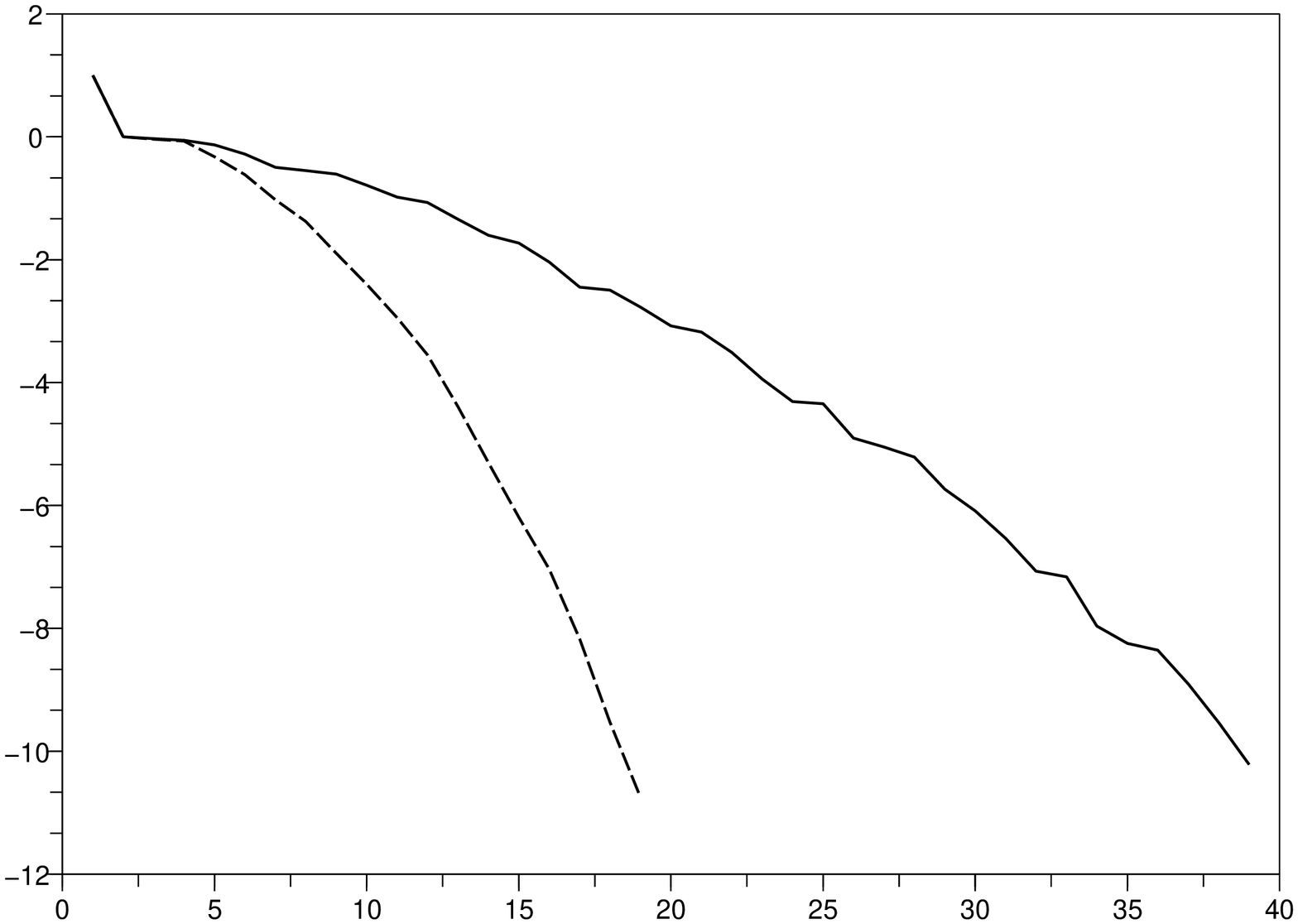}\\Five repetitions
   \end{minipage} \hfill
   \begin{minipage}[c]{.47\linewidth}\centering
      \includegraphics[angle=-90,width=0.99\textwidth]{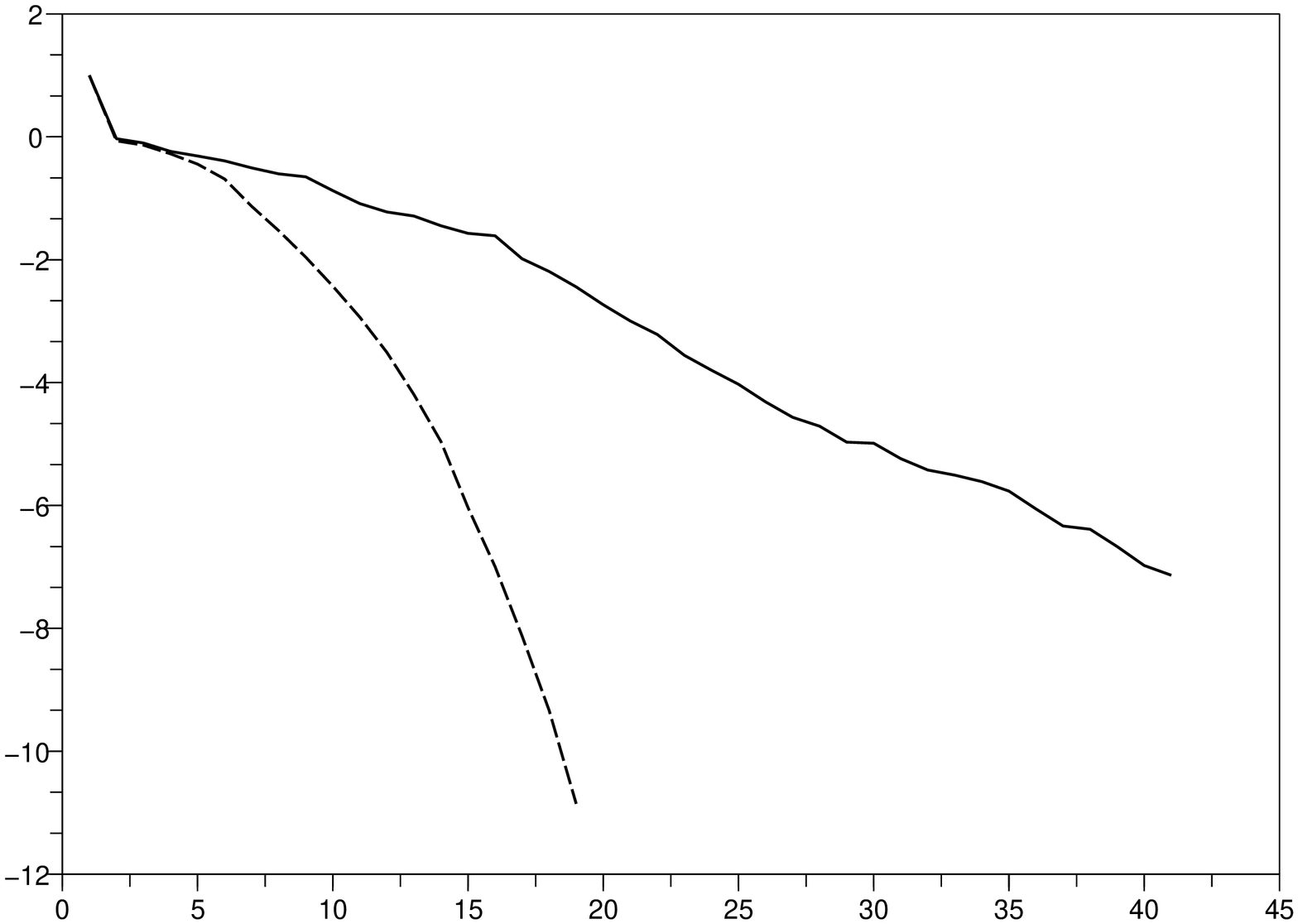}\\Nine repetitions
   \end{minipage}\caption{Comparison between new (dashed) and classical (solid) algorithms}\label{fig:aca_result:1}
\end{figure}

In the present work the method is assessed on academical matrices (randomly generated as in the examples of section
\ref{subsubsec_aca_per}). Each interface is of dimension 20 and
either 5 or 9 repetitions are considered. In figure \ref{fig:aca_result:1} we plot the convergence history of the
classical and the multivector FETI.
It is observed that the exploitation of the repetition leads to high efficiency, the number of iterations being divided by a
factor greater than 2. It is also remarkable that the multivector method is independent on the number of repetitions.

It thus seems that the proposed multivector FETI approach leads to convergence efficiency that is scalable with the
number of repetitions.
A possible interpretation of this result is that the exploitation of the repetition has the same algorithmic effect
as the use of a coarse grid problem (typically associated to rigid body motions of floating substructures):
at each iteration the residual is made orthogonal to a subspace augmented by the permutations
(that is to say a space that is larger than the single vector of classical conjugate gradient).
The permutations enable to transmit instantaneously a global information on the whole structure
which is precisely the aim of coarse grid problems.

\subsubsection{Donut with two stands}
In the case of two identical stands (figure \ref{fig:donut:n}.b) new partial periodicity appears.
It is possible to consider the following enriched multivector:
\begin{equation}
w =
\begin{pmatrix}
w^{(1,2)} \\
w^{(2,3)} \\
w^{(3,1)} \\
w^{(1,5)} \\
w^{(2,4)} \\
\end{pmatrix}
\qquad
\muv{w} =
\begin{pmatrix}
w^{(1,2)} & w^{(2,3)} & w^{(3,1)} & w^{(1,2)} & w^{(2,3)} & w^{(3,1)}\\
w^{(2,3)} & w^{(3,1)} & w^{(1,2)} & w^{(2,3)} & w^{(3,1)} & w^{(1,2)}\\
w^{(3,1)} & w^{(1,2)} & w^{(2,3)} & w^{(3,1)} & w^{(1,2)} & w^{(2,3)}\\
w^{(1,5)} & w^{(1,5)} & w^{(1,5)} & w^{(2,4)} & w^{(2,4)} & w^{(2,4)}\\
w^{(2,4)} & w^{(2,4)} & w^{(2,4)} & w^{(1,5)} & w^{(1,5)} & w^{(1,5)}
\end{pmatrix}
\end{equation}

Of course one danger of this philosophy is an exponential increase of components in the multivector
as the number of patterns and repeated interfaces grows. Proper strategies to choose the permutations
are yet to be developed in that case.

\subsubsection{Dirichlet bc's}
Let us now consider the case of one piece of the donut having part of its outer frontier submitted to Dirichlet
conditions (imposed temperature, figure \ref{fig:donut:n}.c); in that case the Schur complement of the subdomain
should be modified. A first option would be to consider that the structure is made out of two patterns
(the two free pieces versus the clamped one). Another option, following the \textit{Total FETI} method proposed
in \cite{DOSTAL.2006.1}, is to use Lagrange multipliers to impose the boundary conditions and thus the Schur complements
are always associated to free substructures so that all subdomain can be represented by the same pattern.

\subsubsection{The preconditioning problem}
In any of the last two cases we showed that non-periodic interface conditions (extended to the boundary of the subdomain)
could be integrated into our acceleration framework. One key point which gives precedence to the dual approach
is that the dual Schur complement involves the resolution of Neumann problems associated to matrix ${\stiff\s}^+$
which is independent on the location of the imposed fluxes, so that a block treatment of the
$\dgschur \muv{w}$ product is always possible.

The problem is different when dealing with the Dirichlet preconditioner since it is strongly
dependent on the location of the Dirichlet conditions (which modifies the splitting between
interior and boundary degrees of freedom). Thus the theoretically optimal preconditioner is made out
of contributions which are potentially different even for subdomains associated to the same pattern. As an example,
consider again the structure with one stand (figure \ref{fig:donut:n}.a):
the Dirichlet operator of subdomain $3$ is different from the one of subdomains $1$ and $2$
because it should integrate the action of the stand. Thus a block treatment of the optimal preconditioner is complex to manage;
anyhow lumped strategies are still available for block treatment and should lead to good performance results.

\subsection{Linear elasticity problems}
For now we have considered thermal problems because the research of a scalar unknown leads to simpler writing.
To use previous strategy for elastic problems one has to take into account the rotations which are necessary
to transform the pattern to one of its occurrences. These rotations  modify the interpretation of the degrees of freedom
of the interface nodes. Nevertheless this additional local rotation step does not modify the strategy proposed here
(see for instance \cite{RIXEN.2005.1} for a detailed formulation of periodic elastic problems in the dual setting).

More precisely two modifications are required:
\begin{itemize}
 \item The assembly operator needs to incorporate rotation from local to global frame. Since Schur complements are defined in the pattern's frame, not only the assembly operator selects degrees of freedom associated to the interfaces of specific occurence but it expresses them in the global frame: ${\assem\s}^T$ goes from the structure's frame to the pattern's frame and $\assem\s$ sends back data defined on the pattern to data defined on the structure.
 \item The permutations which enable to obtain multivectors ($w$ to $\muv{w}$) also need to take into account the rotation. Data defined on interface $\partial^{1,2}$ in the structure's frame are not directly pertinent for interface $\partial^{2,3}$, they need to be transformed to be defined in the same frame.
\end{itemize}
These operations are easy to implement and do barely incur additional computational costs.

As an example we consider the ``fully-periodic'' donut as a bidimensional elasticity problem (plane strain). The internal ring is clamped and loading is the result of a random process. The pattern has about $2 000$ degrees of freedom ($1 000$ nodes) of which about 100 are on its interface. The acceleration due to the use of multivectors is very similar to the one observed for thermal problems, performance are summed up in table \ref{tab:donut2}. Here again the use of multivectors seems to provide numerical scalability, CPU gain settles around $25 \%$ (which could be improved if some redundant operations where suppressed in the yet non-optimal implementation of the multivector approach).

\begin{table}[ht]\centering
 \begin{tabular}{|c|p{3.1cm}|p{3cm}|}\hline
\multicolumn{3}{|c|}{5-part donut}\\\hline
&  classical FETI-mrhs & multivector FETI   \\\hline
number of iterations &13&10\\  \hline
CPU time (s) &45.25&34.82\\\hline\hline
\multicolumn{3}{|c|}{9-part donut}\\\hline
 & classical FETI-mrhs & multivector FETI \\\hline
number of iterations &15&10\\  \hline
CPU time (s) &43.23&30.09\\\hline
 \end{tabular}\caption{Performance of the multivector approach for elastic donut problems}\label{tab:donut2}
\end{table}


\subsection{Floating substructures}
Part of the treatment of rigid body motions is simplified in our framework because it can be realized at the scale
of the pattern (and not at the level of  its occurrences). For instance in the fully-periodic elastic donut problem, assuming the pattern is only fixed by a hinge at the central node of its internal ring, each occurence has one rigid body motion (rotation around the hinge) and the structure has enough Dirichlet boundary conditions (see figure \ref{fig:pat_rotul}).

\parbox[ht]{0.2\textwidth}{\centering
\includegraphics[width=\linewidth]{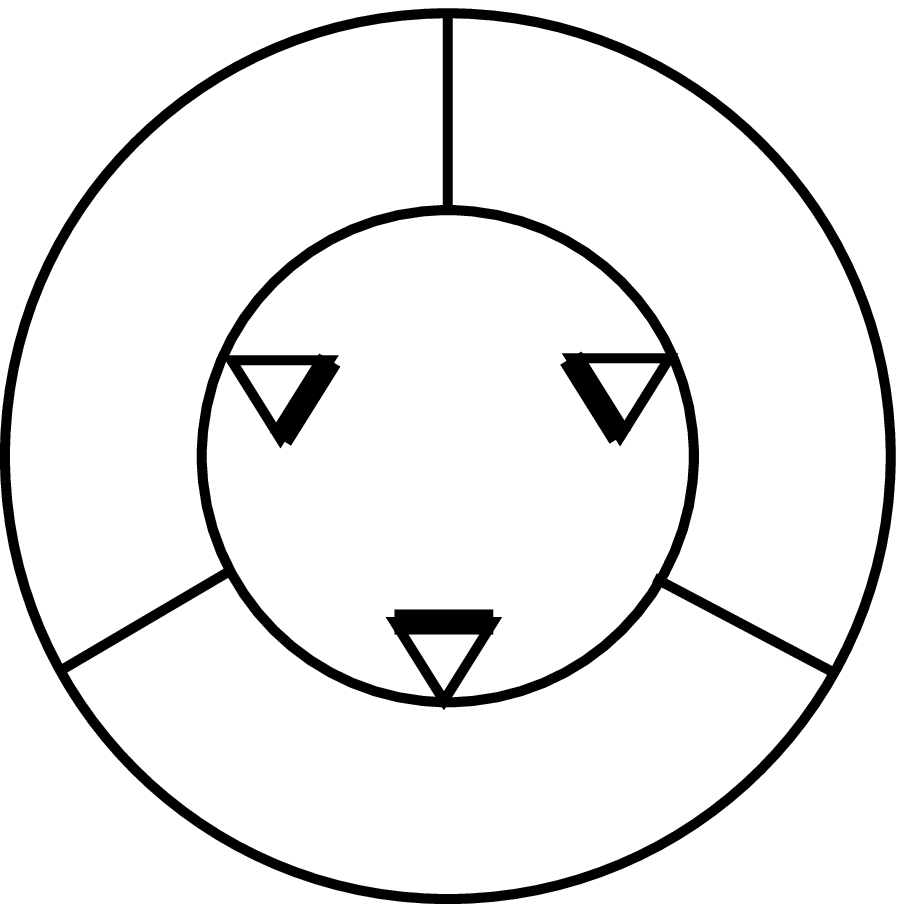}
\captionof{figure}{3-piece structure with hinges}\label{fig:pat_rotul}
} \hfill
\parbox[ht]{0.72\textwidth}{\centering
\begin{tabular}{|c|p{3.1cm}|p{3.cm}|}\hline
\multicolumn{3}{|c|}{5-part donut}\\\hline
&  classical FETI-mrhs & multivector FETI   \\\hline
nb. iterations &12&10\\  \hline
CPU time (s) &44.72&36.23\\\hline\hline
\multicolumn{3}{|c|}{9-part donut}\\\hline
 & classical FETI-mrhs & multivector FETI \\\hline
nb. iterations &14&10\\  \hline
CPU time (s) &38.27&28.28\\\hline
 \end{tabular}
\captionof{table}{Elastic donut problems with floating substructures - }\label{tab:donut3}}

Let $\trace\n\kernel\n=\begin{pmatrix} \kernel_a \\ \kernel_b \end{pmatrix}$ where $a$ and $b$ stand for the two interfaces of the pattern, the matrix $G$ of the trace of the rigid body modes in FETI writes (see equation
(\ref{eq:DUALequilibre})):
\begin{equation}
 G=\begin{pmatrix}
    G^{(1,2)}  & G^{(2,1)} & 0 \\
    0          &  G^{(2,3)} & G^{(3,2)} \\
    G^{(1,3)} & 0          &  G^{(3,1)}
   \end{pmatrix} = \begin{pmatrix}
   \mathcal{O}_a^{(1)} \kernel_a  & -\mathcal{O}_b^{(2)}\kernel_b & 0 \\
    0          &  \mathcal{O}_a^{(2)}\kernel_a & -\mathcal{O}_b^{(3)}\kernel_b \\
    -\mathcal{O}_b^{(1)}\kernel_b & 0          &  \mathcal{O}_a^{(3)}\kernel_a
   \end{pmatrix}
\end{equation}
where matrix $\mathcal{O}_a^{(1)}$ is the rotation matrix which transforms interface $a$ of the pattern into interface $a$ of occurence $(1)$. 

The computation of the coarse grid operator $(G^T\tilde{\dgschur}^{-1}G)$ sums up to the evaluation of
the reaction of the pattern to traces of rigid body modes on its interfaces. So one just has to introduce matrix
$G_{\text{block}}$ defined in the frame of the pattern:
\begin{equation}
G_{\text{block}}=\begin{pmatrix} \kernel_a & \kernel_a & 0 \\ \kernel_b & 0 & \kernel_b \end{pmatrix}
\end{equation}
and deduce $(G^T\tilde{\dgschur}^{-1}G)$ from the computation of
$G_{\text{block}}^T (\scal \tilde{S}\n \scal^T) G_{\text{block}}$ ($\scal$ is the block adaptation of the scaling matrix $\scal\n$, $\tilde{S}\n$ is the chosen approximation of the primal Schur complement of the pattern). Hence the coarse grid problem can be built by computing the following responses of the pattern: its response when submitted to a rigid body mode on its entire interface and its response to the trace of the rigid modes on each of its interface segment successively.

Results presented in table \ref{tab:donut3} (obtained with $\tilde{S}\n=I$) show that rigid body motions have no influences on the performance of the multivector technique. Such a result is in good agreement with previous analysis: rigid body motions are associated to a large wavelength information whereas permutations deal with small wavelength phenomena so that the two effects reinforce one another.

\section{Conclusion and outlook}
This paper presents the first developments around a new approach to compute structures made out of repeated patterns.
The solution algorithm is based on a nonoverlapping domain decomposition method (in order to extract the patterns
from the mesh) and on block Krylov solvers (in order to distribute time consuming operations and
optimize the solution procedure on a larger subspace). The method always converges at least as fast as classical
method with a cost per iteration only slightly higher. We have illustrated the efficiency of
the proposed strategy on simple examples.

Few drawbacks need to be investigated yet, mostly the problem of the generation of redundant information
(which leads to handling too much data and to numerical difficulties).

The extension of the method to more realistic problems is a rather complex task which at first implies an
object oriented description of the mesh. The entire system is then based on the patterns and
connections between their different instantiations.
This requires some adaptation of usual codes and probably a better cooperation with CAD modelling tools (where the information related to the existence of the patterns and their occurences is directly available). Moreover, it is important to note that the same way rotations are used to transport information in the correct frames, non-boolean assembly operators (as obtained with mortar techniques) could be employed to deal with non-matching meshes and transport information from one discretization to another, which would avoid to impose too complex contraints on the meshing process (patterns could be meshed independently from one another, and faces of the same pattern would not need to have the same discretization).

Nevertheless we believe that handling such a description would lead to high computational benefits: only patterns would be meshed (saving memory) and would communicate with one another larger amount of data at one time (saving message passing time),
operations would always be realized on blocks and minimization would be realized on larger subspaces.

\section*{Acknowledgements}
P. Gosselet wishes to thank Erasmus European Exchange Program for partly funding his visit to Delft University of Technology where part of this work was realized.

\bibliographystyle{plain}
\bibliography{ddpattern2}
\end{document}